\documentclass[conference, compsoc]{IEEEtran}

\ifCLASSINFOpdf
\else
  \usepackage[dvips]{graphicx}
\usepackage{psfig}
\fi
\def\RR{I\mskip-7muR}  
\def\ee{I\mskip-7muE} 
\def\X{\mathcal{X}}

\def\bx{\mbox{\boldmath $x$}}  
\def\bX{\mbox{\boldmath $X$}}
\def\bu{\mbox{\boldmath $u$}}

\def\bp{\mbox{\boldmath $p$}}
\def\bk{\mbox{\boldmath $k$}}

\def\btheta{\mbox{\boldmath $\theta$}}
\def\bbeta{\mbox{\boldmath $\beta$}}
\def\bSigma{\mbox{\boldmath $\Sigma$}}

\begin{document}
%
\title{Numerical studies of the metamodel fitting and validation processes}


\author{\IEEEauthorblockN{Bertrand Iooss}
\IEEEauthorblockA{Electricit\'e de France - EDF R\&D\\
Chatou, France\\
bertrand.iooss@edf.fr}
\and
\IEEEauthorblockN{Lo\"{\i}c Boussouf}
\IEEEauthorblockA{Altran\\
Toulouse, France\\
loic.boussouf@gmail.com}
\and
\IEEEauthorblockN{Vincent Feuillard}
\IEEEauthorblockA{EADS IW\\
Suresnes, France\\
vincent.feuillard@eads.net}
\and
\IEEEauthorblockN{Amandine Marrel}
\IEEEauthorblockA{Institut Fran\c{c}ais du P\'etrole\\
Rueil-malmaison, France\\
amandine.marrel@ifp.fr}
}

\maketitle

\begin{abstract}
Complex computer codes, for instance simulating physical phenomena, are often too time expensive to be directly used to perform uncertainty, sensitivity, optimization and robustness analyses.
A widely accepted method to circumvent this problem consists in replacing cpu time expensive computer models by cpu inexpensive mathematical functions, called metamodels. 
In this paper, we focus on the Gaussian process metamodel and two essential steps of its definition phase.
First, the initial design of the computer code input variables (which allows to fit the metamodel) has to provide adequate space filling properties.
We adopt a numerical approach to compare the performance of different types of space filling designs, in the class of the optimal Latin hypercube samples, in terms of the predictivity of the subsequent fitted metamodel.
We conclude that such samples with minimal wrap-around discrepancy are particularly well-suited for the Gaussian process metamodel fitting.
Second, the metamodel validation process consists in evaluating the metamodel predictivity with respect to the initial computer code.
We propose and test an algorithm, which optimizes the distance between the validation points and the metamodel learning points in order to estimate the true metamodel predictivity with a minimum number of validation points.
Comparisons with classical validation algorithms and application to a nuclear safety computer code show the relevance of this new sequential validation design.
\end{abstract}

{\bf Keywords -}
Metamodel, Gaussian process, discrepancy, optimal design, Latin hypercube sampling, computer experiment.

\IEEEpeerreviewmaketitle

\section{Introduction}

With the advent of computing technology and numerical methods, investigation of computer code experiments remains an important challenge.
Complex computer models calculate several output values (scalars or functions), which can depend on a high number of input parameters and physical variables.
These computer models are used to make simulations as well as predictions, uncertainty analyses or sensitivity studies \cite{derdev08}.

However, complex computer codes are often too time expensive to be directly used to conduct uncertainty propagation studies or global sensitivity analysis based on Monte Carlo methods.
To avoid the problem of huge calculation time, it can be useful to replace the complex computer code by a mathematical approximation, called a metamodel \cite{sacwel89,klesar00}.
Several metamodels are classically used: polynomials, splines, generalized linear models, or learning statistical models like neural networks, regression trees, support vector machines \cite{fanli06}.
One particular class of metamodels, the Gaussian process (Gp) model, extends the kriging principles of geostatistics to computer experiments by considering the correlation between two responses of a computer code depending on the distance between input variables \cite{sacwel89}.
Numerous studies have shown  that this interpolating model provides a powerful statistical framework to compute an efficient predictor of code response \cite{sanwil03,marioo07}.

From a practical standpoint, fitting a Gp model implies estimation of several hyperparameters involved in the covariance function. This optimization problem is particularly difficult in the case of a large number of inputs \cite{fanli06,marioo07}.
Several authors (for example \cite{simpep01} and \cite{fanli06}) have shown that the space filling designs are well suited to metamodel fitting.
However, this class of design, which aims at obtaining the better coverage of the points in the space of the input variables, is particularly large, ranging from the well known Latin Hypercube Samples to low discrepancy sequences \cite{fanli06}.
At the moment, no theoretical result gives the type of initial design, which leads to the best fitted Gp metamodel in terms of metamodel predictivity.
In this work, we propose to give some numerical results in order to answer to this fundamental question.

Another important issue we propose to address concerns the optimal choice of the test sample, i.e., the set of simulation design, which allows the most accurate metamodel validation using the minimal number of additional test observations.
The validation of a metamodel is an essential step in practice \cite{klesar00}.
By estimating the metamodel predictivity, we obtain a confidence degree associated with the use of the metamodel instead of the initial numerical model.
Two validation methods are ordinarily used: the test sample approach \cite{ioovan06} and the cross validation method \cite{mecboo02,reipor06}.
In this paper, we propose to perform numerical studies of the metamodel predictivity with respect to these validation methods.

In the following section, we present the Gp model.
In the third section, we present several criteria to optimize the choice of the initial input design.
On two analytical examples, we evaluate the numerical performance of this optimal design in terms of Gp metamodel predictivity.
In the fourth section, we look at the metamodel validation problem.
Our solution consists in minimizing the number of test observations by using the recent algorithm of \cite{feu07}, called the sequential validation design.
We illustrate the relevance of this new design by performing intensive simulation on two analytical functions and an industrial example.
Finally, a conclusion summarizes our results and gives some perspectives for this work.

\section{Gaussian process metamodeling}

Let us consider $n$ realizations of a computer code. Each realization $y(\bx) \in \RR$ of the computer code output corresponds to a $d$-dimensional input vector $\bx = (x_1,\ldots,x_d) \in \X$, where $\X$ is a bounded domain of $\RR^d$. The $n$ points corresponding to the code runs are called the experimental design and are denoted as $\bX_s = ( \bx^{(1)},\ldots,\bx^{(n)} )$. The outputs will be denoted as $Y_s = (y^{(1)},\ldots,y^{(n)})$ with $y^{(i)} = y(\bx^{(i)})$ $\forall \; i=1..n$.
Gaussian process (Gp) modeling treats the deterministic response $y(\bx)$ as a realization of a random function $Y(\bx)$, including a regression part and a centered stochastic process. This model can be written as:
\begin{equation}\label{Gpmodel}
 Y ( \bx ) = f ( \bx ) + Z ( \bx) .
\end{equation}

The deterministic function $f(\bx)$ provides the mean approximation of the computer code. In our study, we use a one-degree polynomial model where $f(\bx)$ can be written as follows:
\[ f ( \bx ) = \beta_0 + \sum_{j = 1}^d \beta_j x_j \;, \]
where $\bbeta = [ \beta_0, \ldots, \beta_k ]^t $ is the regression parameter vector.
It has been shown, for example in \cite{marsim05} and \cite{marioo07}, that such a function is sufficient, and sometimes necessary, to capture the global trend of the computer code.

The stochastic part $Z(\bx)$ is a Gaussian centered process fully characterized by its covariance function:
$\mbox{Cov} ( Z ( \bx ), Z ( \bu ) ) = \sigma^2 R ( \bx, \bu ),$
where $\sigma^2$ denotes the variance of $Z$ and $R$ is the correlation function that provides interpolation and spatial correlation properties.
To simplify, a stationary process $Z(\bx)$ is considered, which means that correlation between $Z(\bx)$ and $Z(\bu)$ is a function of the distance between $\bx$ and $\bu$. Our study is focused on a particular family of correlation functions that can be written as a product of one-dimensional correlation functions $R_l$:
\[ \mbox{Cov} ( Z ( \bx ), Z ( \bu ) )  = \sigma^2 R ( \bx - \bu ) = \sigma^2 \prod_{l = 1}^d R_l ( x_l - u_l ) . \]
This form of correlation functions is particularly well adapted to get some simplifications of integrals in analytical uncertainty and sensitivity analyses \cite{marioo09}.
More precisely, we choose to use the generalized exponential correlation function:
\[ R_{\btheta, \bp} ( \bx - \bu ) = \prod_{l = 1}^d \exp ( - \theta_l |x_l - u_l |^{p_l} )  ,\]
where $\btheta = [ \theta_1, \ldots, \theta_d ]^t $ and $\bp =  [ p_1, \ldots, p_d ]^t $ are the correlation parameters (also called hyperparameters) with $\theta_l \geq 0$ and $0 < p_l \leq 2$ $\;\forall\; l=1..d$.
This choice is motivated by the wide spectrum of shapes that such a function offers.

If a new point $\bx^{\ast} = (x^{\ast}_1,\ldots,x^{\ast}_d) \in \X$ is considered, we obtain the predictor and variance formulas:
 \begin{eqnarray}
 \displaystyle \ee [  Y_{\mbox{\scriptsize Gp}}(\bx^{*})] = f(\bx^{*}) +  \bk(\bx^{*})  ^t \bSigma_s^{-1} (Y_s - f(\bX_s)) \;, \label{eq_esperance} \\ 
 \displaystyle \mbox{Var}[ Y_{\mbox{\scriptsize Gp}}(\bx^{*})] = \sigma^2  -  \bk(\bx^{*}) ^t  \bSigma_s^{-1} \bk(\bx^{*}) \;, \label{eq_variance}
\end{eqnarray}
with $Y_{\mbox{\scriptsize Gp}}$ denoting $(Y|Y_s,\bX_s,\bbeta,\sigma,\btheta, \bp)$,
  \[ 
\begin{array}{lll}
\bk(\bx^{*}) & = &[\mbox{Cov}(y^{(1)},Y(\bx^{*})), \ldots, \mbox{Cov}(y^{(n)},Y(\bx^{*}))  ] ^t  \\
& = & \sigma^2  [  R_{\btheta, \bp} (\bx^{(1)},\bx^{*}), \ldots, R_{\btheta, \bp} (\bx^{(n)},\bx^{*}) )  ]^t  
 \end{array} 
\]
and the covariance matrix
\[ \bSigma_s = \sigma^2 \left( R_{\btheta, \bp} \left( \bx^{(i)} - \bx^{(j)} \right)_{i=1..n, j = 1..n}  \right) \;.\]
The conditional mean (Eq. (\ref{eq_esperance})) is  used as a predictor. The variance formula (Eq. (\ref{eq_variance}) corresponds to the mean squared error (MSE) of this predictor and is also known as the kriging variance. This analytical formula for MSE gives a local indicator of the prediction accuracy. More generally, Gp model provides an analytical formula for the distribution of the output variable at any arbitrary new point. This distribution formula can be used for sensitivity and uncertainty analysis \cite{marioo09}. 

Regression and correlation parameters $\bbeta$, $\sigma$, $\btheta$ and $\bp$ are ordinarily estimated by maximizing likelihood functions \cite{fanli06}.
This optimization problem can be badly conditioned and difficult to solve in high dimensional cases ($d>5$) \cite{marioo07}.
Moreover, the estimation algorithms are particularly sensitive to the input design.
The following section proposes to deal with this input design problem.

\section{Initial design for the metamodel fitting}\label{sec:fitting}

For computer experiments, selecting an experimental design is a key issue in building an efficient and informative metamodel.
In this section, we describe the different properties than a computer experimental design has to reach.
Some numerical tests support our discussion.

\subsection{Latin hypercube sampling}

\begin{figure}[!ht]
$$\psfig{figure=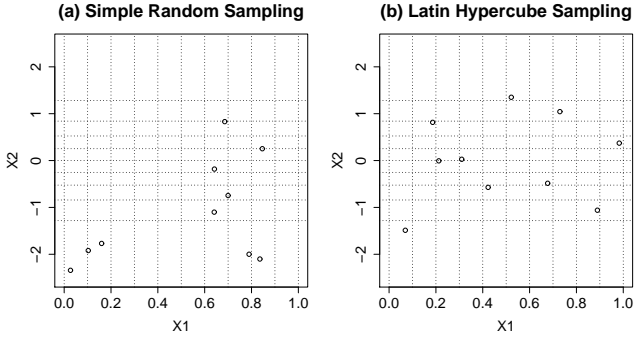,height=6.cm,angle=-90}$$

\vspace{-1.2cm}
\caption{Examples of two ways to generate a sample of size $n=10$ from two variables $\mathbf{X} = [X_1,X_2]$ where $X_1$ has a uniform distribution ${\cal U}[0,1]$ and $X_2$ has a normal distribution ${\cal N}(0,1)$. Equprobable stratas are shown in each dimension.}\label{fig:comparaisonSRSLHS}
\end{figure}

Contrary to the Simple Random Sample (SRS, also called crude Monte Carlo sample), which consists of $n$ independently and identically distributed samples, the well known Latin Hypercube Sample (LHS) consists in dividing the domain of each input variable in $n$ equiprobable strata, and in sampling once from each stratum \cite{mckbec79}.
The LHS of a random vector $\bX=(X_1,\ldots,X_d)$, denoted $(\bX^{(1)},\ldots,\bX^{(n)})$, gives a sample mean $m=\frac{1}{n}\sum_{i=1}^n Y^{(i)}$ for the output $Y=y(\bX)$ with a smaller variance than the sample mean of a SRS \cite{ste87}.
Figure ~\ref{fig:comparaisonSRSLHS} shows $10$ samples of two random variables, $X_1$ and $X_2$, obtained with SRS and LHS schemes.
We can see that the result of LHS is more spread out and does not display the clustering effects found in SRS.

However, LHS does not reach the smallest possible variance for the sample mean.
Since it is only a form of stratified random sampling and it is not directly related to any criterion, it may also perform poorly in metamodel estimation and prediction of the model output.
Therefore, some authors have proposed to enhance LHS not only to fill space in one dimensional projection, but also in higher dimensions \cite{par93}.
One powerful idea is to adopt some optimality criterion applied to LHS, such as entropy, integrated mean square error, minimax and maximin distances, etc.
For instance, the maximin criterion consists in maximizing the minimal distance between the points \cite{johmoo90}.
This leads to avoid situations with too close points.
The paper \cite{mormit95} examines some optimal maximin distance designs constructed within the class of Latin hypercube arrangements.
The conceptual simplicity of these designs has led to their large popularity in practical applications \cite{jonjoh09}.

\subsection{Low-discrepancy Latin hypercube samples}

Alternative metamodel-independent criteria, based on discrepancy measures, consist in judging the uniformity quality of the design.
Discrepancy can be seen as a measure between an initial configuration and an uniform one. 
It is a comparison between the volume of intervals and the number of points within these intervals \cite{hic98}.
There exists different kinds of definition using different forms of intervals or different norms in the functional space.
Discrepancy measures based on $L_2$ norms are the most popular in practice because they can be analytically expressed and are easy to compute.
Among them, two measures have shown remarkable properties \cite{jinche05,fan01,fanli06}:
 \begin{itemize}
 \item
the centered $L^2$ discrepancy
{\scriptsize
\begin{equation}\label{disccent}
\hspace{-1cm}
\begin{array}{l}
\displaystyle D^2(\bX_s(n)) =  \left(\frac{13}{12}\right)^d -\frac{2}{n}\sum_{i=1}^n\prod_{k=1}^d\left(1+\frac{1}{2}|u_k^{(i)}-\frac{1}{2}|-\frac{1}{2}|u_k^{(i)}-\frac{1}{2}|^2\right)\\
\displaystyle +\frac{1}{n^2}\sum_{i,j=1}^{n}\prod_{k=1}^{d}\left(1+\frac{1}{2}|u_k^{(i)}-\frac{1}{2}|+\frac{1}{2}|u_k^{(j)}-\frac{1}{2}|-\frac{1}{2}|u_k^{(i)}-u_k^{(j)}|\right)
\end{array}
\end{equation}}
\noindent where $\bX_s(n)$ denotes the input learning sample with $n$ input vectors and $\left(u_k^{(i)}\right)_{i=1..n,k=1..d}$ are the normalized values in $[0,1]$ of the design $\bX_s(n)=\left(x_k^{(i)}\right)_{i=1..n,k=1..d}$;
\item
the wrap-around $L^2$ discrepancy
{\scriptsize 
\begin{equation}\label{discwrap}
\hspace{-1cm}
W^2(\bX_s(n)) = \left(\frac{4}{3}\right)^d +\frac{1}{n^2}\sum_{i,j=1}^n\prod_{k=1}^d\left[\frac{3}{2}-|u_k^{(i)}-u_k^{(j)}|(1-|u_k^{(i)}-u_k^{(j)}|)\right]\;,
\end{equation}}
\noindent which allows to suppress bound effects (by wrapping the unit cube for each coordinate).
\end{itemize}

\begin{figure}[!ht]
$$\psfig{figure=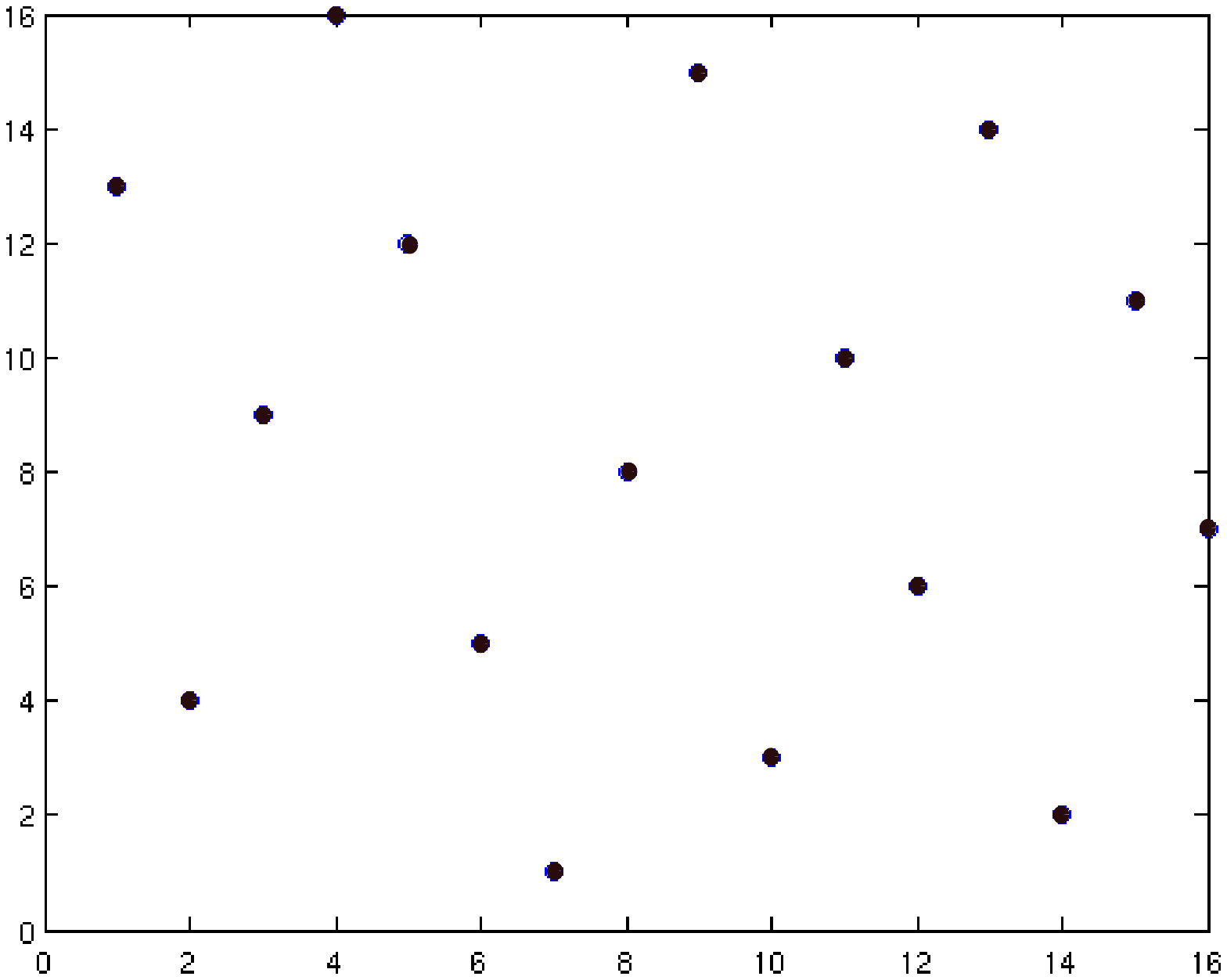,height=3.5cm}
\psfig{figure=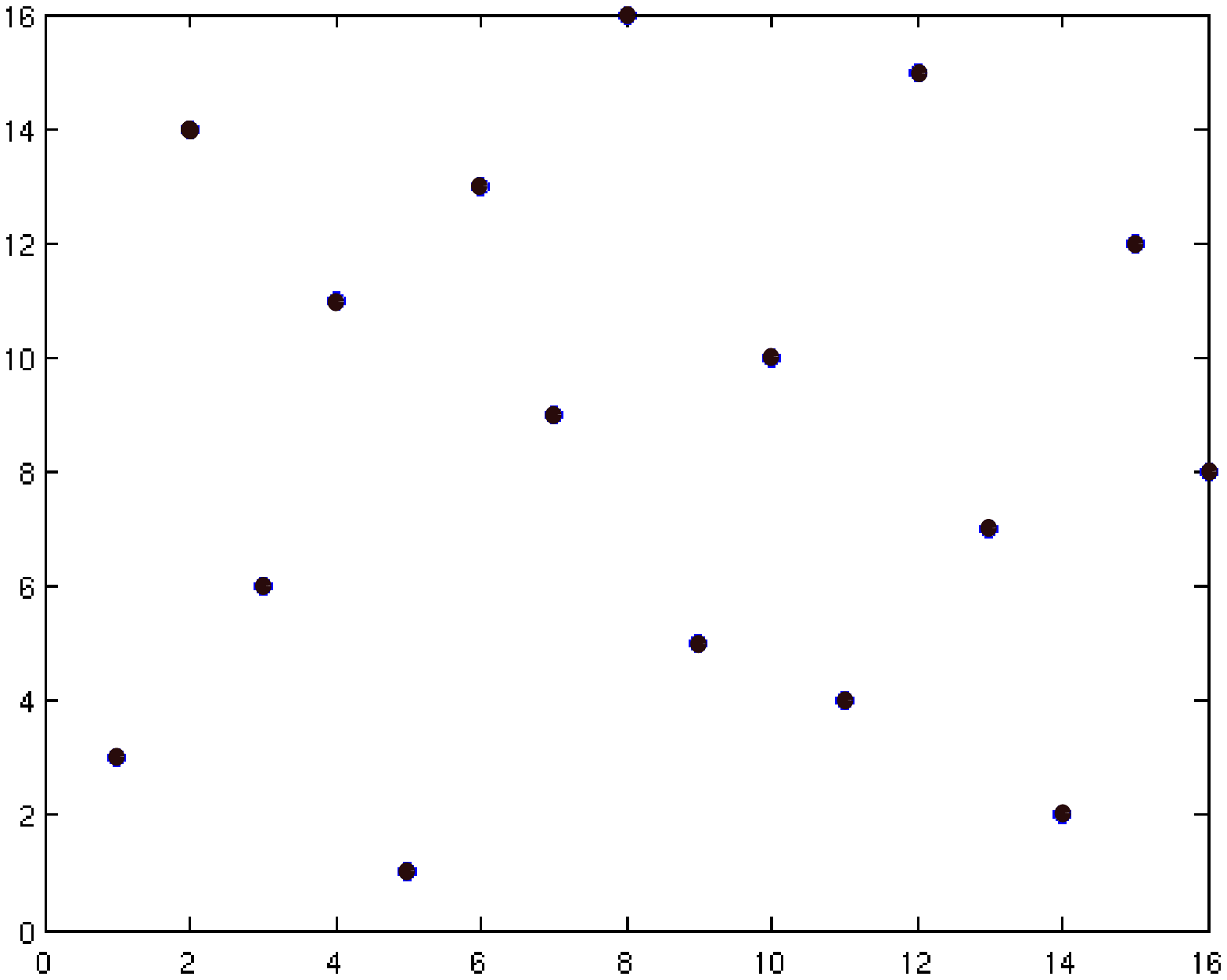,height=3.5cm}$$

\vspace{-0.3cm}
\hspace{1.9cm}Maximin \hspace{1.9cm} Centered discrepancy
$$\psfig{figure=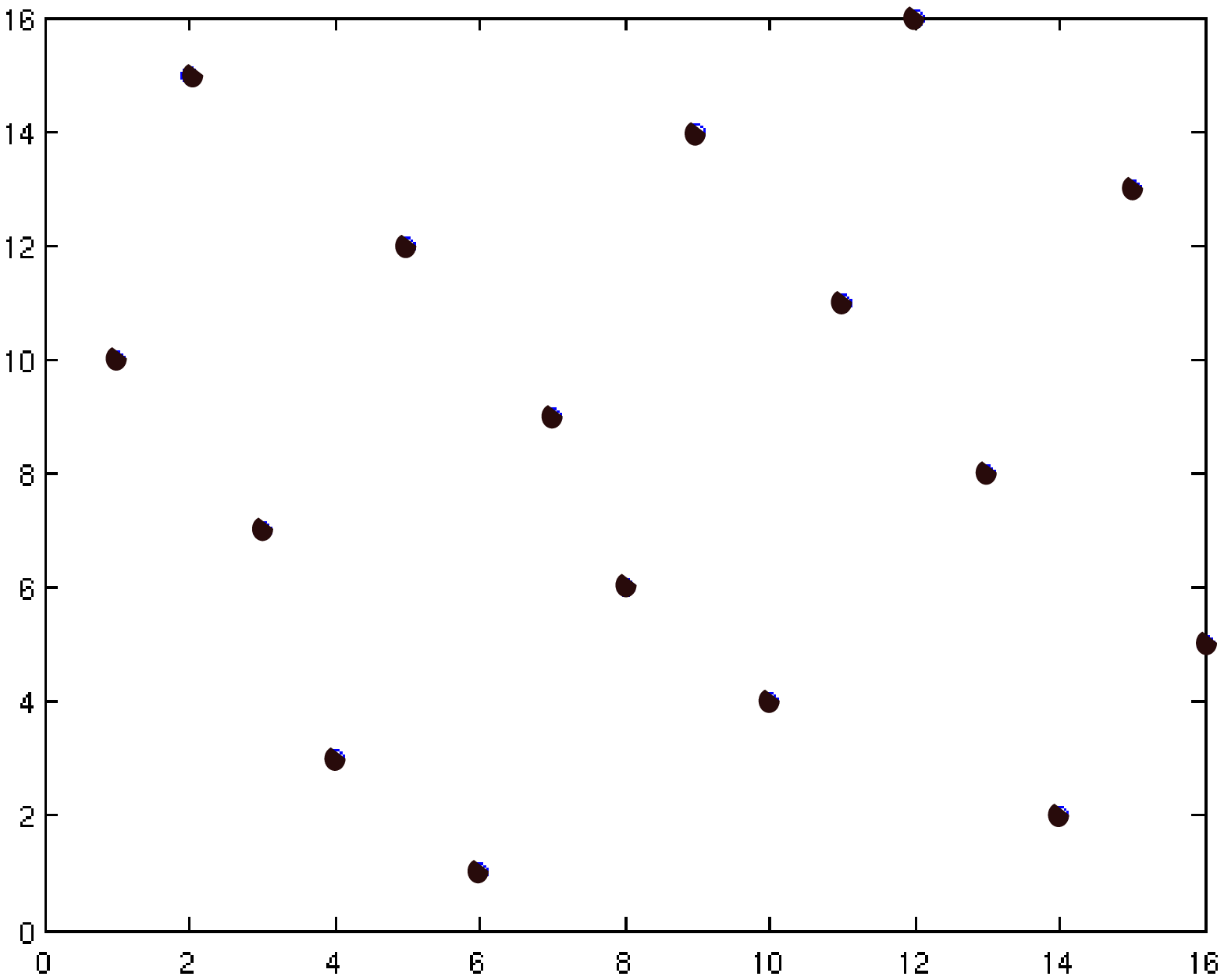,height=3.5cm}$$

\vspace{-0.3cm}
\centerline{Wrap-around discrepancy}
\caption{Visual comparisons of LHS ($d=2$, $n=16$) optimized following three different criteria (below each figure).}\label{LHScriteria}
\end{figure}

The optimization of LHS can be done following different methods: choice of the best (in terms of the chosen criteria) LHS amongst a large number of different LHS, columnwise-pairwise exchange algorithms, genetic algorithms, simulated annealing, etc \cite{jinche05,liesto06}.
In our tests, we have found that the simulated annealing algorithm (with a geometrical temperature descent and with a slight noise on the initial condition) gives the best results for all the criteria \cite{mar08}.
Figure \ref{LHScriteria} gives some examples of two-dimensional LHS of size $n=16$, optimized following three different criteria with the simulated annealing algorithm.
We see that uniform repartitions of the points are nicely respected.

\subsection{Projection properties of space filling designs}

In addition to the space filling property on the sample space, one important property of the initial designs is their robustness to the dimension decrease.
A LHS structure for the space filling design is not sufficient because it only guarantees good repartitions for one-dimensional projections, and not for the other dimensions of projection.
Indeed, LHS ensures that each of the input variables has all proportion of its range which is represented (equiprobable stratas are created for each input variable).
In contrary, no equiprobable stratas are created in the various multi-dimensional spaces of the input variables.

We then argue that the sample points of a space filling design have to be well spread out when projected onto a subspace spanned by a subset of coordinate axes.
This property is particularly important when the initial design is made in dimension $d$ and the metamodel fitting is made in a smaller dimension (see an example in \cite{cangar08}).
In practice, this is often the case because the initial design may reveal with screening methods the useless (i.e., non influent) input variables that we can neglect during the metamodel fitting step \cite{puj09}.
Moreover, when a selection of input variables is made during the metamodel fitting step (as for example in \cite{marioo07}), the new sample, solely including the retained input variables, has to keep good space filling properties.

\begin{figure}[!ht]
$$\psfig{figure=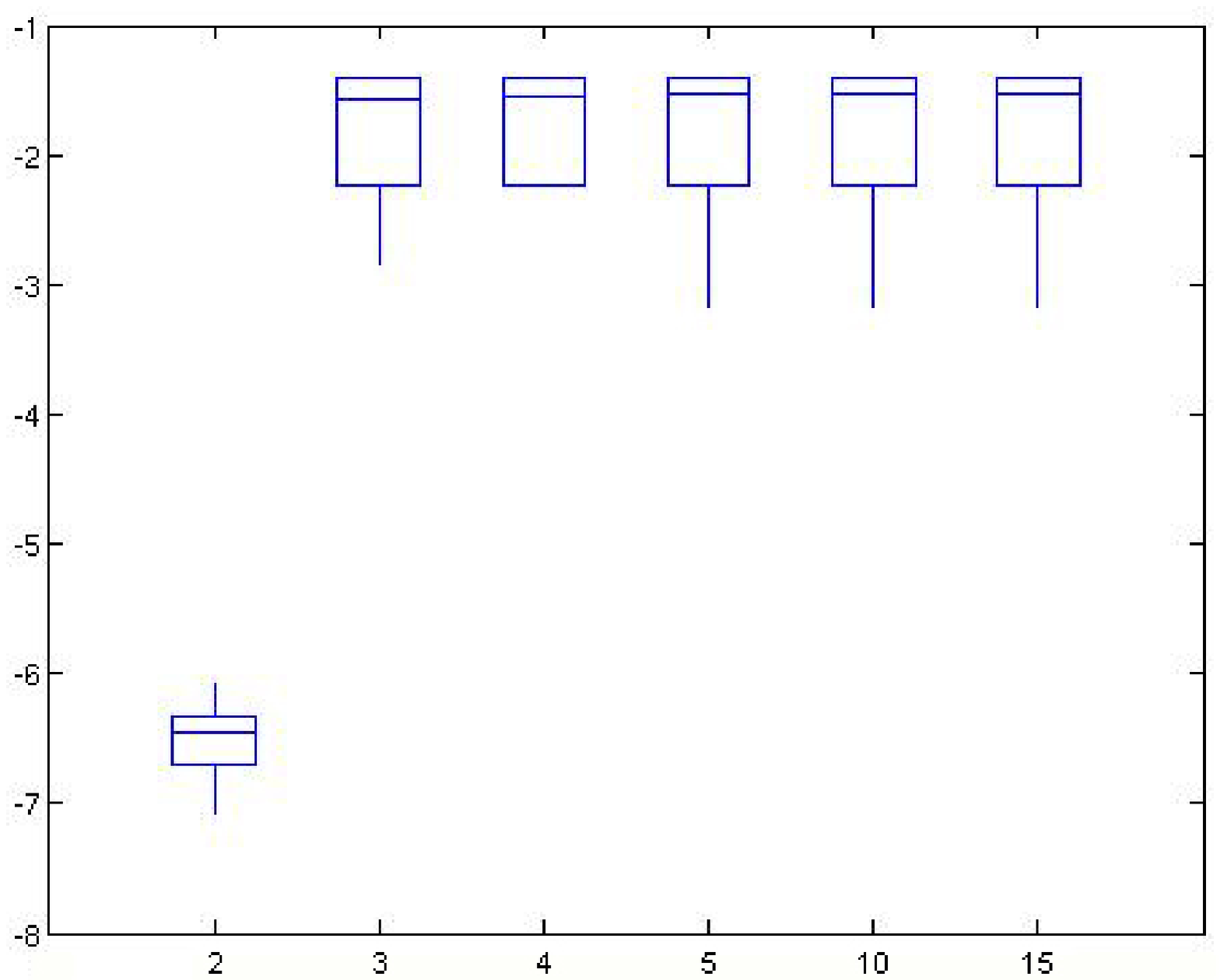,height=4cm,width=7cm,angle=-90}$$

\vspace{-0.3cm}
 \centerline{Maximin LHS}
$$\psfig{figure=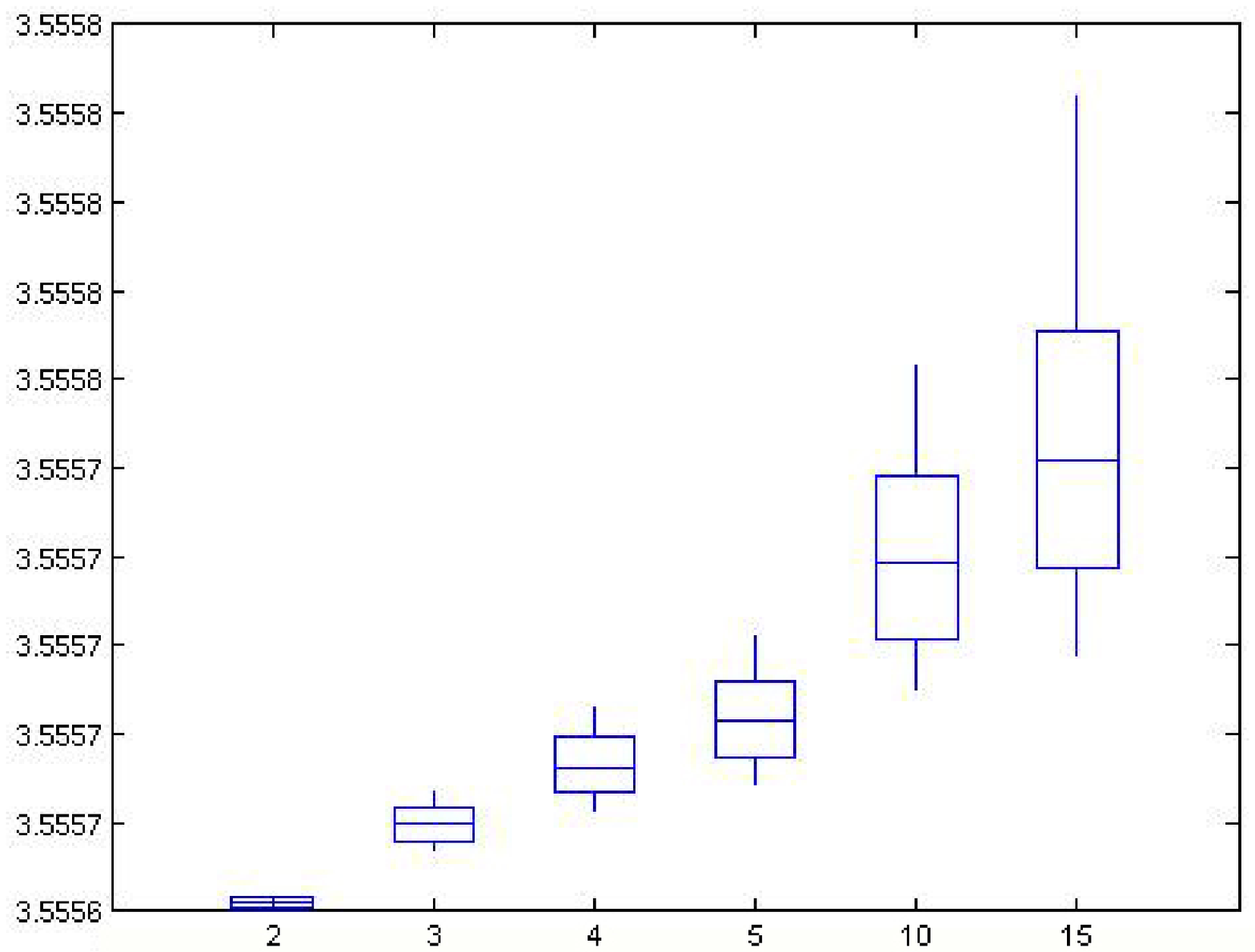,height=4cm,width=7cm,angle=-90}$$

\vspace{-0.3cm}
\centerline{Low wrap-around discrepancy LHS (WLHS)}
\caption{Criterion values (up: maximin, bottom: wrap-around discrepancy) obtained with 2D projections of designs coming from two types of LHS (containing $n=100$ points), with different dimensions: $d=2, 3, 4, 5, 10, 15$. Boxplots are obtained by repeating $100$ optimizations using different initial LHS.}\label{LHScritproj2d}
\end{figure}

Figure \ref{LHScritproj2d} compares the two-dimensional projections of the maximin LHS and low wrap-around discrepancy LHS (called WLHS) with $n=100$ points and different initial dimensions (from $d=3$ to $15$).
The reference criterion values are given for $d=2$.
For dimension larger than $2$, we compute the new criterion values by considering all the two-dimensional projections of the initial design.
A robust criterion to the dimension decrease would lead to a small increase of the criterion value.
The criteria behave very diferently between the two types of design:
\begin{itemize}
\item 2D projection criteria of WLHS regularly and slightly deteriorate. 
Then, 2D projections of WLHS made in dimensions close to $2$ keep rather good space-filling properties.
\item 2D projection criteria of maximin LHS sharply and strongly deteriorate from the first dimension increase at $d=3$.
Then 2D projection criteria of maximin LHS remain stable at poor values for larger dimensions.
\end{itemize}
Similar tests with different sample sizes $n$ and for the three-dimensional and four-dimensional projections have led to the same conclusions.
All these results show that a WLHS is the preferable initial design for fitting a computer code metamodel in high dimensional cases.

\subsection{Numerical studies on toy functions}

At present, we perform two numerical studies to evaluate the impact of an inadequate design on the metamodel fitting process.
For the metamodel, we use the Gp model $Y_{\mbox{\scriptsize Gp}}$ described in \S 2.
The quality of the metamodel predictor is measured by the so-called predictivity coefficient $Q_2$ (i.e., the determination coefficient $R^2$ computed on a test sample), which gives the percentage of the output variance explained by the metamodel:
\begin{equation}\label{Q2}
Q_2 = 1 - \frac{\sum_{i=1}^{n_t} [y(\tilde{\bx}^{(i)})-\hat{Y}_{\mbox{\scriptsize Gp}}(\tilde{\bx}^{(i)})]^2}{\sum_{i=1}^{n_t} [\bar{y}-y(\tilde{\bx}^{(i)})]^2}
\end{equation}
with $(\tilde{\bx}^{(1)},\ldots,\tilde{\bx}^{(n_t)})$ the test sample of size $n_t$, $\hat{Y}_{\mbox{\scriptsize Gp}}=\ee(Y_{\mbox{\scriptsize Gp}})$ the Gp predictor (Eq. (\ref{eq_esperance})) and $\bar{y}$ the mean of the output test sample $(y(\tilde{\bx}^{(1)}),\ldots,y(\tilde{\bx}^{(1)}))$.

\subsubsection{A two-dimensional test case}

\begin{figure}[!htb]  
$$\psfig{figure=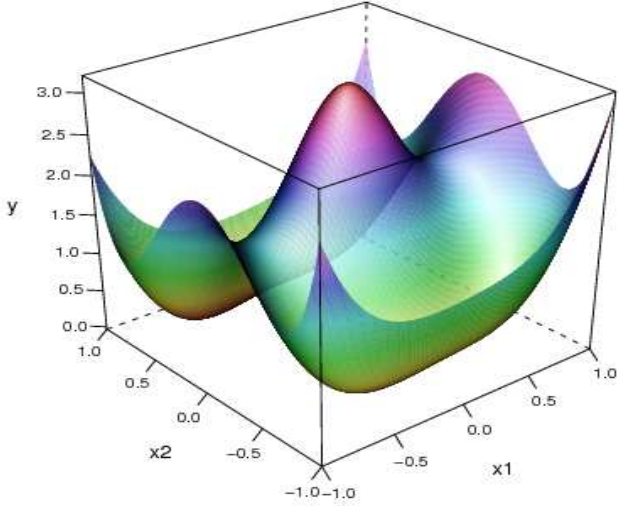,height=7cm,width=9cm}$$

\vspace{-0.5cm}
\caption{Graphical representation of the irregular function on $[-1;1]^2$.}
\label{fig10}
\end{figure}

Our first test involves a two-dimensional analytical function (called the irregular function):
$$ f(\bx) = \frac{e^{x_1}}{5}-\frac{x_2}{5}+\frac{x_2^6}{3}+4x_2^4-4x_2^2+\frac{7x_1^2}{10}+x_1^4+\frac{3}{4x_1^2+4x_2^2+1} $$
with $\bx \in [-1 , 1]^2$.
Figure \ref{fig10} represents the irregular function.

We have made several comparisons between random LHS and different space filling designs before fitting a metamodel \cite{mar08}.
In the following, we show our results concerning the random LHS and the WLHS, which has provided the best results.
For a size $n$ of the learning sample and each type of design, we repeat $100$ times the following procedure: we generate an initial input design of $n$ observations, we obtain $n$ outputs with the toy function, we fit a Gp metamodel (\ref{Gpmodel}), and we evaluate its predictivity coefficient $Q_2$ using a test sample of large size ($n_t=10000$).
Therefore, for each type of LHS, we obtain $100$ values of $Q_2$ whose mean and variance give us the efficiency and robustness of the design in terms of Gp quality.

\begin{figure}[!htb]  
\begin{center}
\includegraphics[height=3.7cm, keepaspectratio=true]{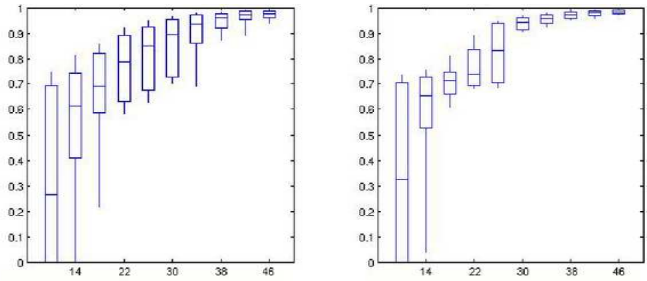}
\end{center}
\caption{For the irregular function, Gp $Q_2$ evolution in function of the learning sample size $n$ and for two types of LHS (left: random LHS; right: WLHS).}
\label{ameliorationQ2}
\end{figure}

The initial LHS design optimized with the wrap-around discrepancy (Eq. (\ref{discwrap})) has given us the best results.
In Figure \ref{ameliorationQ2}, we compare the predictivity coefficients obtained with non optimized LHS (random LHS) and those obtained with optimized LHS (WLHS).
The size of the design increases from $n=10$ to $n=46$ (by step of $4$), which leads to a regular increase of $Q_2$.
For each size $n$, the boxplot represents the summary of the $100$ values of $Q_2$.
In the all range of $n$, $Q_2$ of the WLHS are better than the random LHS ones.
Furthermore, much smaller variances (boxplots are smaller) are shown for WLHS and lead to the conclusion that these designs are more robust than others.
This property is rather natural because there are much less variability between the $100$ different WLHS than between the $100$ different random LHS (because of the optimization process).
Differences are particularly important for sizes $n=30$ and $n=34$: the WHS lead to very competitive Gp metamodels ($Q_2 \sim 0.95$ and boxplot width $\sim 0.05$) while random LHS give uncompleted metamodels ($Q_2 \sim 0.9$ and boxplot width $\sim 0.2$).

\subsubsection{A five-dimensional test case}

Our second test involves a five-dimensional analytical function (called the g-Sobol 5d function):
$$ f(\bx) = \sum_{i=1}^5 \frac{|4x_i-2|+a_i}{1+a_i} $$
with $a_1=1$, $a_2=2$, $a_3=3$, $a_4=4$, $a_5=5$, $\bx \in [0 , 1]^5$. 

We have made several comparisons between random LHS and different space filling designs before fitting a metamodel \cite{mar08}.
In the following, we show our results concerning the random LHS and the WLHS, which has provided the best results.
For a size $n$ of the learning sample and each type of design, we repeat $100$ times the following procedure: we generate an initial input design of $n$ observations, we obtain $n$ outputs with the toy function, we fit a Gp metamodel (\ref{Gpmodel}), and we evaluate its predictivity coefficient $Q_2$ using a test sample of large size ($n_t=10000$).
Therefore, for each type of LHS, we obtain $100$ values of $Q_2$ whose mean and variance give us the efficiency and robustness of the design in terms of Gp quality.

\begin{figure}[!ht]
$$\psfig{figure=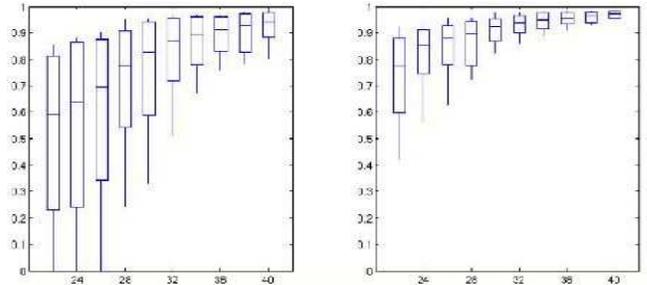,height=3.8cm,width=8.5cm}$$
\caption{For the g-Sobol 5d function, Gp $Q_2$ evolution in function of the learning sample size $n$ and for two types of LHS (left: random LHS; right: WLHS).}\label{LHScritQ2}
\end{figure}

As in the previous section, the initial LHS design optimized with the wrap-around discrepancy (Eq. (\ref{discwrap})) has given us the best results.
In Figure \ref{LHScritQ2}, we compare the predictivity coefficients obtained with non optimized LHS (random LHS) and those obtained with optimized LHS (WLHS).
The size of the design increases from $n=22$ to $n=40$ (by step of $2$), which leads to a regular increase of $Q_2$.
For each size $n$, the boxplot represents the summary of the $100$ values of $Q_2$.
In the all range of $n$, $Q_2$ of the WLHS are better than the random LHS ones.
Furthermore, much smaller variances are shown for WLHS and lead to the conclusion that these designs are more robust than others.
For small sample sizes, the $Q_2$ differences reach $0.2$ between the two types of design: $Q_2(\mbox{LHS}) \sim 0.6$ and $Q_2(\mbox{WLHS}) \sim 0.8$.
In industrial applications, such a difference makes the distinction between ``bad'' (unacceptable) metamodels and good ones.
The latter can be used for example for quantitative sensitivity studies.

\subsection{Conclusion of numerical tests}

In conclusion of our numerical study, the LHS optimized with the wrap-around discrepancy has provided efficient results for the Gp metamodel fitting, even in high dimension.
Furthermore, we have found that this design guarantees correct repartitions of the points for all the two-dimensional projections, while other types of LHS (like maximin) have bad repartitions for these projections.
Other types of LHS can also provide good results but less systematically \cite{mar08}.
For instance, \cite{fra07} has studied quasi-Monte Carlo samples (Sobol suites and Halton sequences) and has shown that these sequences are less performant than other space filling designs in terms of the Gp metamodel fitting.

Of course, such designs have to be seen as initial ones.
If possible, in a second step, adaptive designs can improve metamodel predictivity in a very efficient way \cite{mar08}, for instance by choosing new simulation points in poorly predicted areas.

\section{Test sample selection for metamodel validation}

In practical cases, only a small number of simulations can be performed with the computer code in order to fit a metamodel.
Once the metamodel has been built, estimating its predictivity is an important issue.
Indeed, a safe use of this metamodel to answer to uncertainty or sensitivity problems requires a precise estimation of its capabilities.
In this section, we make a discussion on algorithms of predictivity estimation.

\subsection{Classical validation methods}

Let us consider the $d$-dimensional input vector $\bx = (x_1,\ldots,x_d) \in \X$, where $\X$ is a bounded domain of $\RR^d$ and $y(\bx) \in \RR$ is the computer code output. 
We suppose that a metamodel $\hat{Y}(\bx)$ has been fitted using $\left((\bx^{(1)},y(\bx^{(1)})), \ldots, (\bx^{(N)},y(\bx^{(N)}))\right)$, a $N$-size learning sample of computer code experiments.

The test sample approach consists in comparing the metamodel predictions on simulation points not used in the metamodel fitting process.
This gives some prediction residuals (which can be finely analyzed) and global quality measures as the metamodel predictivity coefficient $Q_2$ (Eq. (\ref{Q2})).
Such test points set is called a test sample (or also validation sample or prediction sample).
This method requires new calculations with the computer code and the first question we have to face up is the sufficient number of prediction points to obtain the required accuracy of our global validation measures.
For cpu time expensive code, it can be difficult to provide a sufficient number of test points.
Some convergence visualisation tools of the global validation measures can be used to answer to this first question.

Another important question for the test sample approach is the localization of these test points.
The usual practice is to choose an independent Monte Carlo sample for the test sample.
However, if the sample size is small, the proposed points can be badly localized, for example near learning points or leaving large space domain unsampled.
A fine strategy could be to use, as the test sample, a space filling design (which consists in filling the input variable space $\X$ as uniformly as possible).
Unfortunately, this solution does not avoid the possibility of too strong proximity between learning points and test points.
Such proximity would lead to too optimistic quality measures, and consequently to a biased prectivity estimation.

The second solution to validate a metamodel, the cross validation method, is extremely popular in practice because it avoids new calculations on the computer code.
The cross validation method proposes to divide the initial sample on a learning sample and a test sample.
A metamodel is estimated with the points in the new learning sample and prediction residuals are obtained via the new test sample.
This process is repeated several times by using other divisions of the learning sample.
Finally all the prediction residuals can be used to compute the global predictivity measures.
The leave-one-out procedure is a particular case of the cross validation method where just one observation is left out at each step.

The first drawback of the cross validation method is its cost, which can become large due to many metamodel fitting processes.
Moreover, if the initial design has a specific geometric structure (which aims at optimizing the metamodel fitting), the deletion of points from the learning sample causes the breakdown of the specific design structure while creating the new learning sample.
Indeed, the new learning sample does not have the adequate statistical and geometric properties of the initial design and the metamodel fitting process might fail. 
This could lead to too pessimistic quality measures.

To sum up, the test sample method requires too many new prediction points (to avoid too optimistic validation criteria), while the cross-validation method can provide 
too pessimistic validation criteria.
Therefore, to solve this dilemma, an heuristic new solution has been introduced in \cite{ioobou08,ioo09} and is presented in the next section.

\subsection{A new optimized validation design}

Retaining the test sample method, we limit its main drawback by minimizing the number of necessary points in the test sample.
In this goal, an algorithm allows the specification of new design points decreasing the discrepancy of an initial design \cite{feu07}.
This sequential algorithm gives us at each step the prediction point furthest away from the other points of the design.
The algorithm performs its optimization process in the space $\X$ of the input variables $\bx$.
By choosing the future prediction points in the unfilled zone of the learning sample design, we aim at capturing the right metamodel predictivity using only a small number of additional points.
Note that such ideas have also been proposed in \cite{ren09} for different purposes.

We have not theoretically studied the computational efficiency of this algorithm over the computational efficiency of the traditional methods (introduced in the previous section).
However, our intuition is that mean square error computed by this algorithm avoids the biases, which could be caused by too strong proximities between the test sample points and between test sample points vs. learning sample points.

Let us consider $\bX_f(n_f)=(\bx_f^{(i)})_{i=1..n_f}$ a low discrepancy sequence of $n_f$ points in $[0,1]^d$.
A low discrepancy sequence is a deterministic design constructed to uniformly fill the space with regular patterns.
Among all the low discrepancy sequence, Halton, Hammersley, Faure and Sobol sequences are the most famous \cite{lem09}. 
In the following, we will use the Hammersley sequence which, on a few tests, have shown better properties than the others \cite{feu07}.
The chosen discrepancy measure is the centered $L^2$ discrepancy $D^2(\cdot)$ (Eq. (\ref{disccent})).

To obtain an additional point of the initial $N$-size sample, noticed $\bX_s(N)$, we use the following algorithm:
\begin{enumerate}
\item For $i=1,\ldots,n_{f}$,
\begin{itemize}
			\item[$\bullet$] $\bX_s(N+1)=\{\bx^{(1)},\ldots,\bx^{(N)}\} \cup {\bx_{f}}^{(i)}$;
      \item[$\bullet$] compute ${\rm Dif}_i = D^2(\bX_s(N+1))-D^2(\bX_s(N))$;
\end{itemize}
\item select $i^*$ such that $\displaystyle i^*=\arg\!\!\!\!\!\min_{i=1,\ldots,n_{f}}{\rm Dif}_{i}$;
\item obtain the new point ${\bx_{f}}^{(i^*)}$.
\end{enumerate}
This algorithm is repeated sequentially to obtain $N_{\mbox{\scriptsize test}}$ test points, by updating the initial design and the low discrepancy sequence.
For example, for the second point, we reinitalize the design by the following: $\bX_s(N+1)=\{\bx^{(1)},\ldots,\bx^{(N)}\} \cup {\bx_{f}}^{(i^*)}$ and $\bX_{f}(n_{f}-1)=\{{\bx_{f}}^{(1)},\ldots,{\bx_{f}}^{(n_f)}\} \backslash {\bx_{f}}^{(i^*)}$.

\begin{figure}[!ht]
\centering
\includegraphics[width=2.5in,keepaspectratio=true]{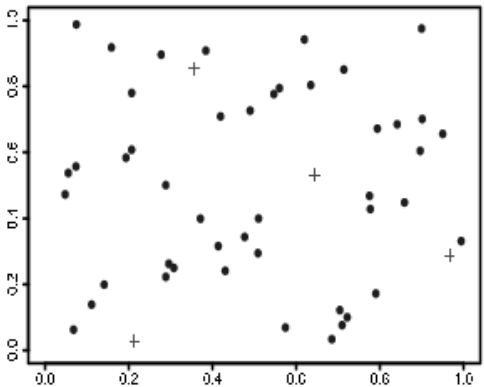}
\caption{Example of the sequential algorithm: $N=46$, $d=2$, $N_{\mbox{\scriptsize test}}=4$. The bullets are the points of the initial design while the crosses are the new specified points.}\label{B1}
\end{figure}

This algorithm just consists in adding to the initial design some points of a low discrepancy sequence by minimizing the discrepancy differences between the initial and the new design.
The size of the low discrepancy sequence is required to be as large as possible, especially if $d$ is large.
Figure \ref{B1} gives an example of the specification with our algorithm of $N_{\mbox{\scriptsize test}}=4$ new points (the crosses) inside an initial Monte Carlo design ($N=46$, $d=2$).
One of the advantage of this algorithm is its size-independence (related to the number of added points): the sequence of added points is deterministic and will be always the same for the same $\bX_f(n_f)$.
In the following, the design obtained using this algorithm is called the sequential validation design.

\subsection{Numerical studies on toy functions}

\subsubsection{A two-dimensional test case}

\begin{figure}[!ht]
$$\psfig{figure=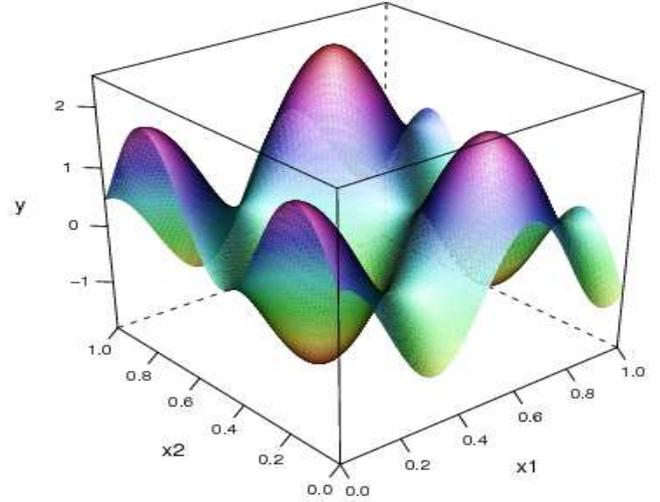,height=7cm,width=9cm}$$

\vspace{-0.5cm}
\caption{Graphical representation of the cosin2 function on $[0;1]^2$.}\label{fig:cosin2}
\end{figure}

To compare the sequential validation design with other test designs for the metamodel validation purpose, we first perform an analytical test using a two-dimensional toy function, called the cosin2 function:
$$f(\bx) = \cos(10x_1)+\sin(10x_2)+x_1 x_2 \;,\; (x_1,x_2) \in [0 , 1]^2.$$
Figure \ref{fig:cosin2} represents the cosin2 function.

Gp metamodels (\ref{Gpmodel}) are fitted using learning samples of differents sizes $N_{\mbox{\tiny BA}}$: $N_{\mbox{\tiny BA}}$ ranges from $10$ to $40$ allowing a wide variety of metamodel predictivity coefficients $Q_2$, from $0$ (null predictivity) to $1$ (perfect predictivity).
The initial $10$-size design is a maximin LHS.
The other learning designs (of increased size) are obtained by sequentially adding points to the design, while maintaining the LHS properties of the design and keeping some optimality properties (maximizing the mean distance from each design point to all the other points in the design \cite{liesto06}).
Choosing an initial maximin LHS design, while we have shown in Section \ref{sec:fitting} that WLHS is better than maximin LHS for the Gp fitting process, is not in contradiction with our objectives in this section: our goal is now to study the Gp metamodel validation.
Anyway, we are not able to keep the properties of maximin LHS or WLHS when gradually increasing the size of the learning sample.

\begin{figure}[!ht]
$$\psfig{figure=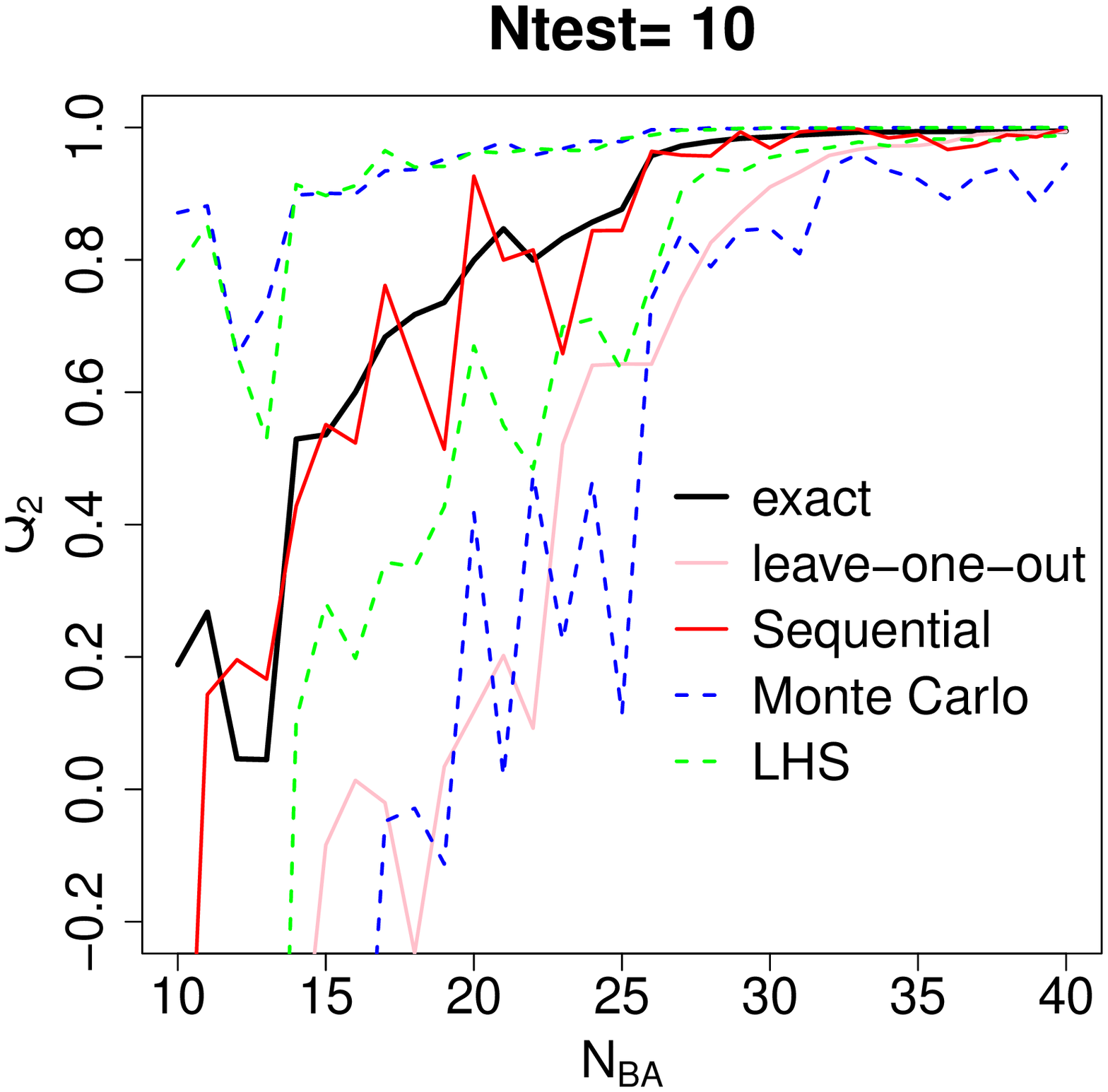,height=4.9cm,width=6.2cm}$$
$$\psfig{figure=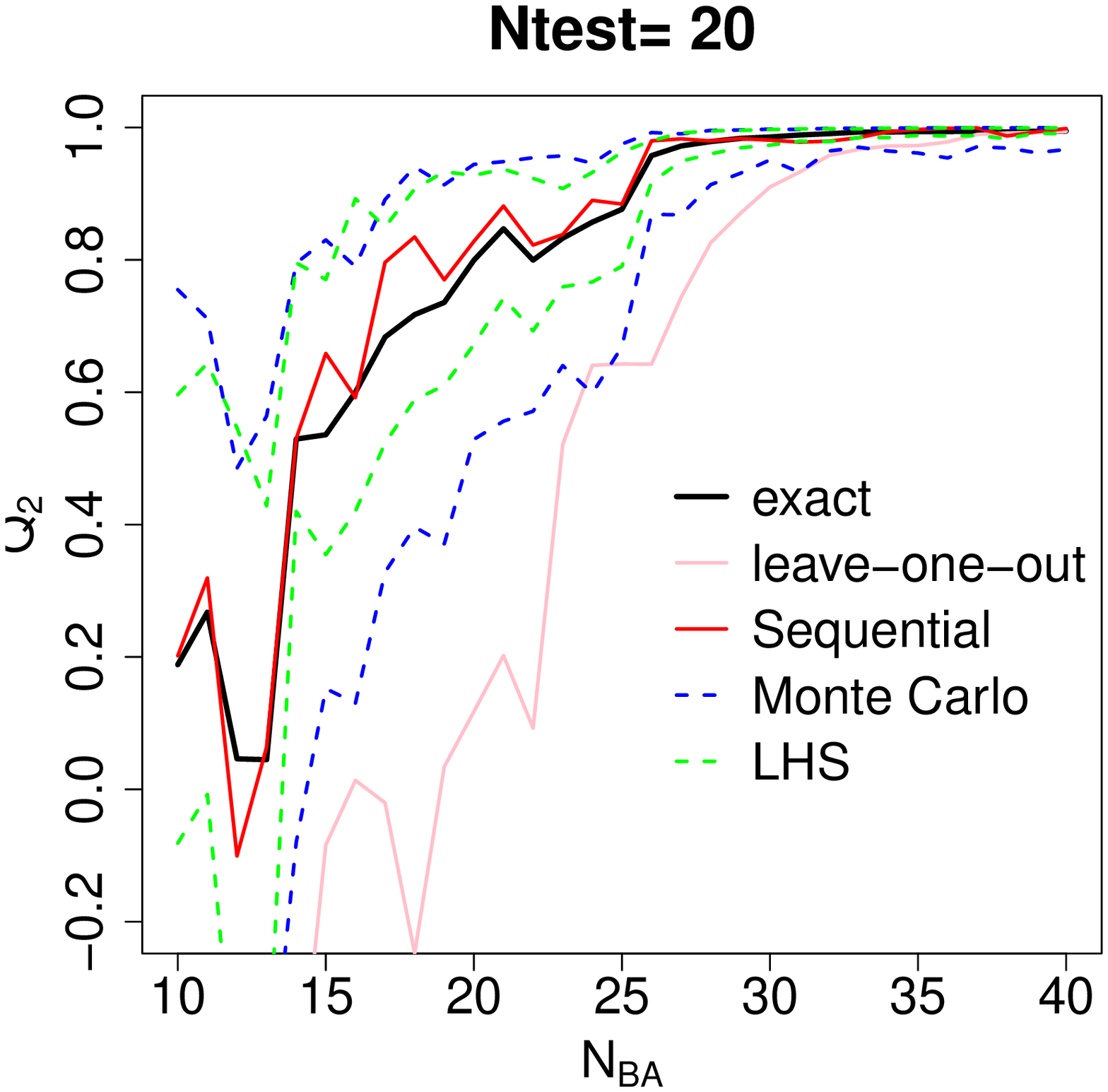,height=4.9cm,width=6.2cm}$$
$$\psfig{figure=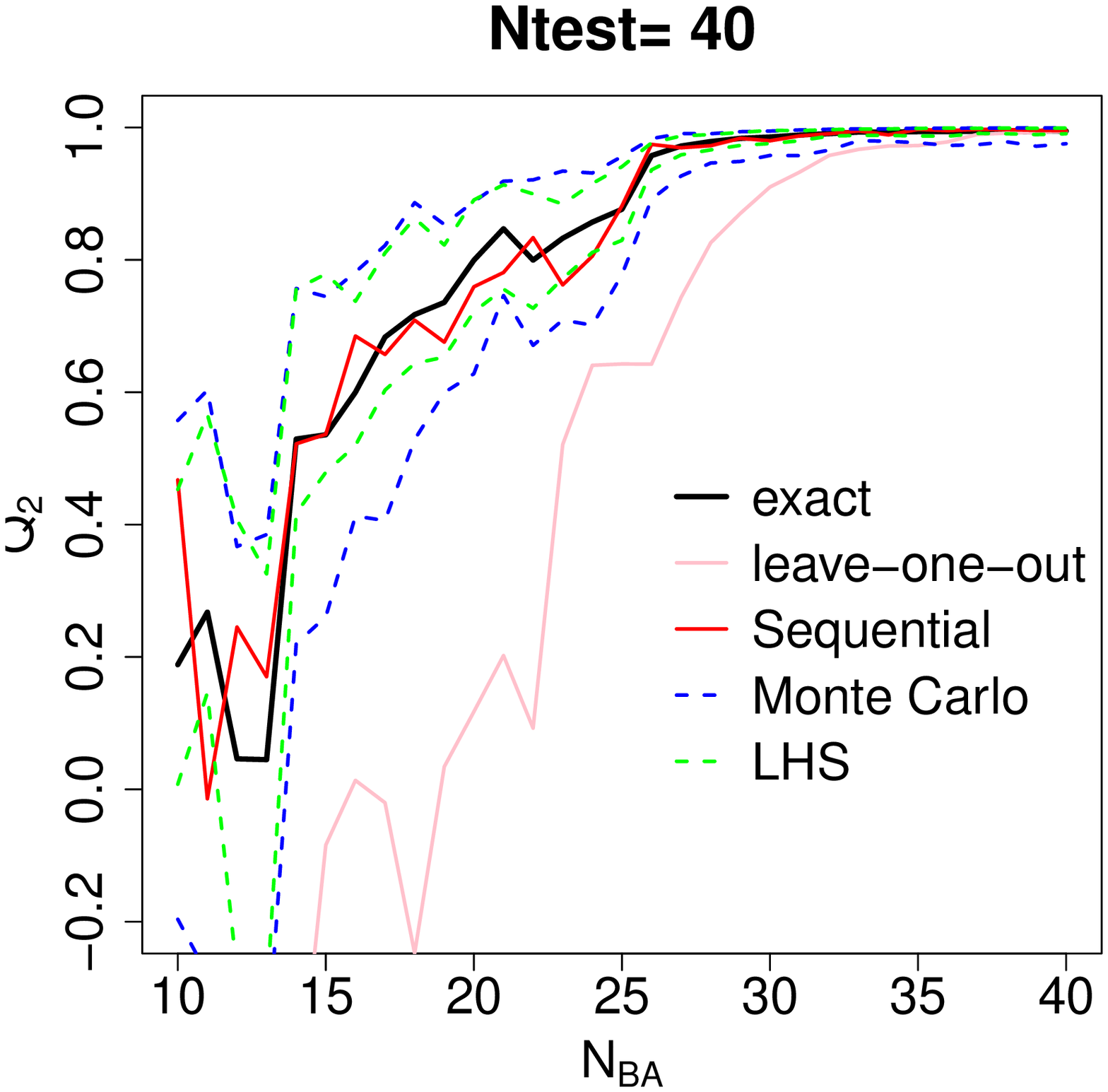,height=4.9cm,width=6.2cm}$$
$$\psfig{figure=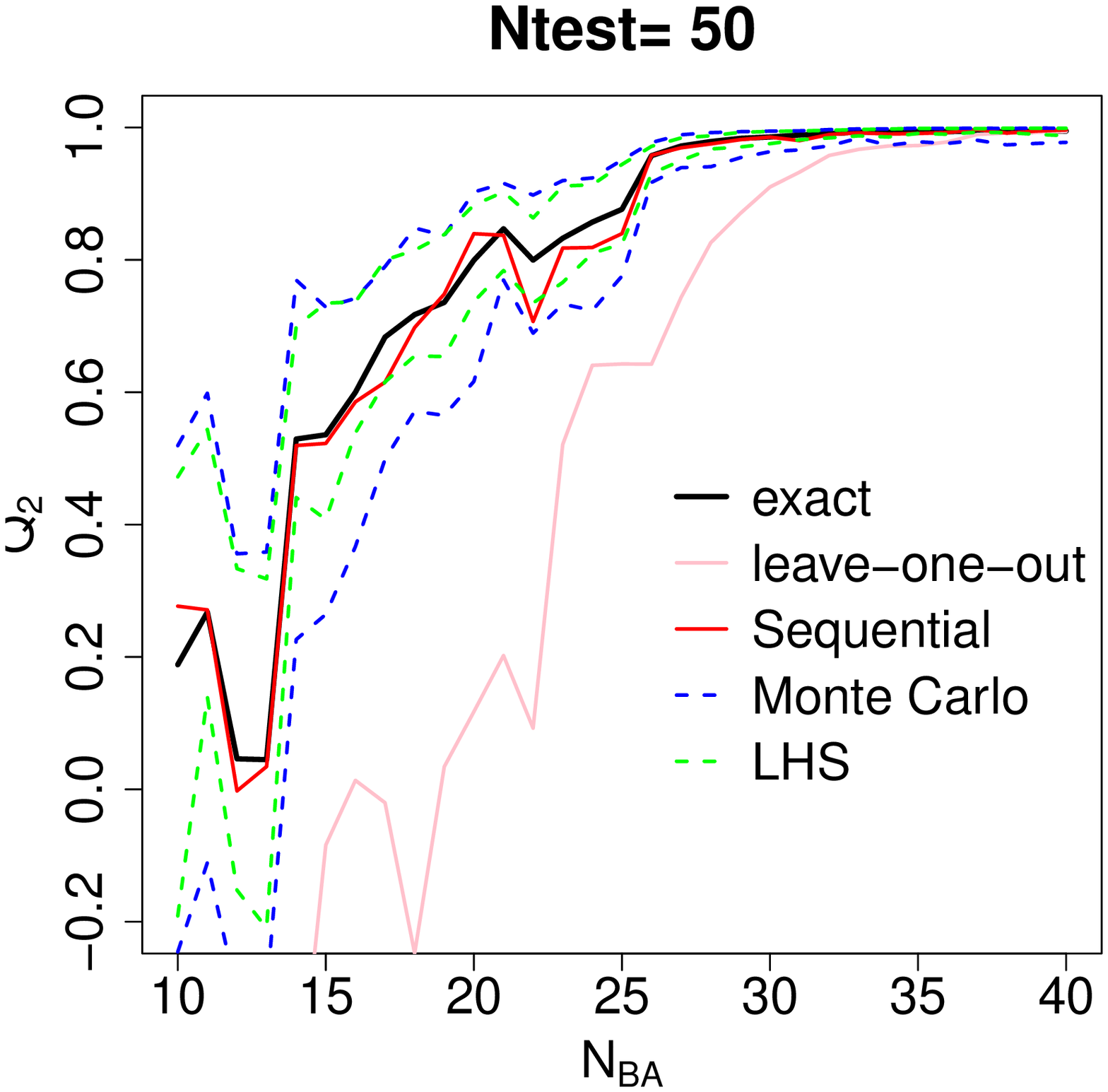,height=4.9cm,width=6.2cm}$$

\vspace{-0.5cm}
\caption{For the cosin2 function, Gp predictivity coefficient ($Q_2$) in function of the learning sample size $N_{\mbox{\tiny BA}}$, estimated from different test sample sizes $N_{\mbox{\scriptsize test}}$. The dashed curves (blue and green) give the minimal and maximal values obtained with $100$ repetitions of the random test design (Monte Carlo and LHS).}\label{Q2toy}
\end{figure}

The black line in Figure \ref{Q2toy} shows the evolution of $Q_2$ in function of the learning sample size.
This reference value for the predictivity coefficient has been computed for each metamodel by taking its mean over $100$ test samples of size $N_{\mbox{\scriptsize test}}=1000$.
The $Q_2$ estimation by a leave one out procedure (pink line) strongly underestimates the exact $Q_2$ for $N_{\mbox{\tiny BA}}<30$.
This is certainly due to the small number of points: leave one out is pessimistic in this case because each point deletion has a strong impact on the metamodel fitting process.
The red curve gives the $Q_2$ estimation using the sequential validation design described in the previous paragraph (with a Hammersley sequence of size $n_f=10000$).

\clearpage

Results are greatly satisfactory for $N_{\mbox{\scriptsize test}} \geq 20$: the sequential validation design gives precise $Q_2$ estimates in all cases and outperforms a crude Monte Carlo or LHS design.
The green curves correspond to the minimal and maximal values obtained with $100$ repetitions using an optimized LHS as the test design.
As expected, these intervals are more reduced than the intervals obtained using a crude Monte Carlo sample as the test design (blue curves).
As $N_{\mbox{\scriptsize test}}$ increases, these intervals contract, but always show the superiority of the sequential validation design, especially for low metamodel predictivity ($Q_2<0.9$ and $N_{\mbox{\tiny BA}}<25$).

\subsubsection{An eight-dimensional test case}

We perform now a second numerical test using the g-Sobol function in eight-dimension (called the g-Sobol 8d function):
$$ f(\bx) = \sum_{i=1}^8 \frac{|4x_i-2|+a_i}{1+a_i}$$
with $a_1=a_2=3$, $a_i=0$ for ($i=3,\ldots,8$), $\bx \in [0 , 1]^8$. 

A Gp metamodel (\ref{Gpmodel}) is fitted on a learning sample (maximin LHS) of size $N_{\mbox{\tiny BA}}=40$.
We compute the reference value of the predictivity coefficient by taking its mean over $100$ test samples of size $N_{\mbox{\scriptsize test}}=1000$ and obtain $Q_2^{\mbox{\scriptsize ref}}=0.83$.
We then apply the sequential validation design described previously (with a Hammersley sequence of size $n_f=10000$) by adding $N_{\mbox{\scriptsize test}}=50$ new points to the design, and we obtain $Q_2^{\mbox{\scriptsize seq$50$}}=0.85$, which is close to the true value.
We compare this result with $100$ crude Monte Carlo samples of the same size ($N_{\mbox{\scriptsize test}}=50$), which give the 90\% confidence interval $[0.79,0.91]$ for $Q_2^{\mbox{\scriptsize MC}}$.
This last result is rather large and shows the insufficient number of points if we choose a crude Monte Carlo design.

\begin{figure}[!ht]
\centering
\includegraphics[width=3in]{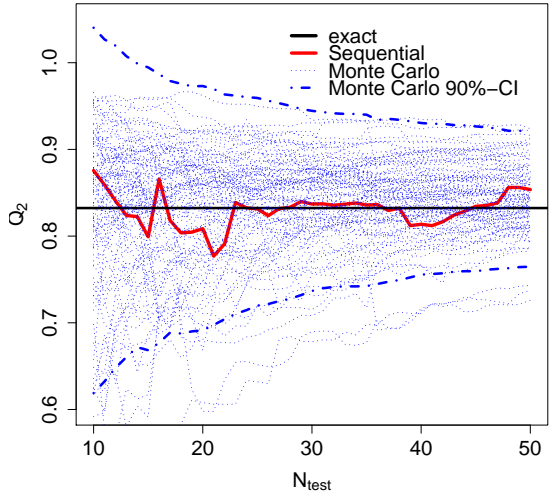}

\vspace{-0.4cm}
\caption{For the g-Sobol 8d function, Gp predictivity coefficient ($Q_2$) in function of the test sample size $N_{\mbox{\scriptsize test}}$, for two types of validation design: sequential (red) and crude Monte Carlo (blue). Dotted blue lines correspond to $100$ different crude Monte Carlo samples).}\label{Q2gsobol}
\end{figure}

Figure \ref{Q2gsobol} shows the evolution of the estimated $Q_2$ for test bases with different sizes, ranging from $N_{\mbox{\scriptsize test}}=10$ to $N_{\mbox{\scriptsize test}}=50$.
The solid red line shows the results obtained with the sequential validation design while the dotted blue lines show the $100$ sequentially increased crude Monte Carlo samples.
This figure illustrates the poor estimates we obtain when using small size ($N_{\mbox{\scriptsize test}}<50$) of Monte Carlo samples for validation.
On the contrary, the sequential validation design allows to obtain a good approximation of the true predictivity coefficient even for small test sample sizes.
Results are precise for $N_{\mbox{\scriptsize test}} \geq 25$.

\subsection{Application to a nuclear safety computer code}

\begin{figure}[!ht]
$$\psfig{figure=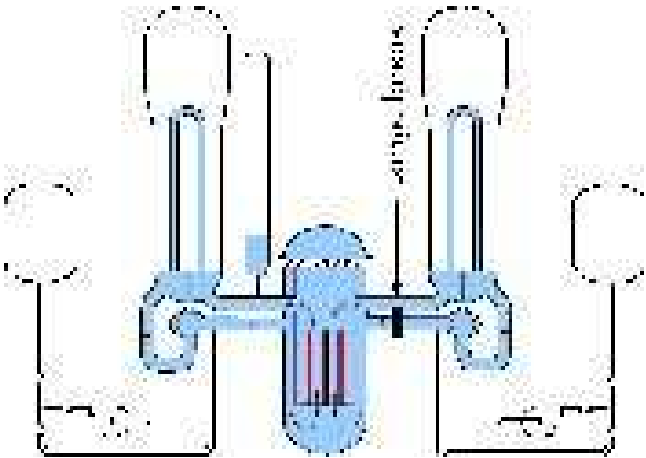,height=6cm,width=8.5cm}$$
\caption{Illustration of a large-break loss of primary coolant accident on a nuclear Pressurized Water Reactor.}\label{bemuse}
\end{figure}

\begin{figure}[!ht]
$$\psfig{figure=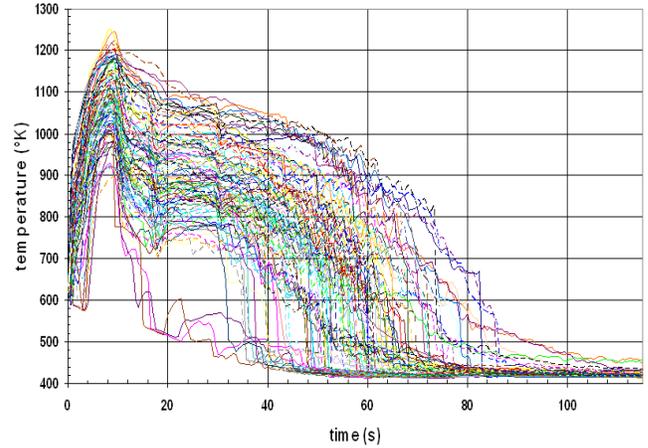,height=6cm,width=8.5cm}$$
\caption{$100$ output curves (cladding temperatures in function of time) of the CATHARE2 code. The output variable of interest for the reactor safety is the first peak of the cladding temperature.}\label{sorties_bemuse}
\end{figure}

In this section we apply our algorithms on a complex computer model used for nuclear reactor safety. It simulates a hypothetical thermal-hydraulic scenario on Pressurized Water Reactors:  a large-break loss of primary coolant accident (see Fig. \ref{bemuse}) for which the output of interest is the peak cladding temperature.
This scenario is part of the Benchmark for Uncertainty Analysis in Best-Estimate Modelling 
for Design, Operation and Safety Analysis of Light Water Reactors \cite{dec08}
proposed by the Nuclear Energy Agency of the Organisation for Economic Co-operation and Development (OCDE/NEA).
It has been implemented on the french computer code CATHARE2 developed at the Commissariat \`a l'Energie Atomique (CEA).
Figure \ref{sorties_bemuse} illustrates $100$ CATHARE2 simulations (by varying input variables of the accidental scenario) giving the cladding temperature in function of time.

In our exercise, a Gp metamodel (\ref{Gpmodel}) of the first peak cladding temperature (which is a scalar variable) has to be estimated with $N=100$ simulations of the computer model (the input design is a maximin LHS). 
The cpu time is twenty minutes for each simulation with a standard computer (Pentium IV PC). 
The complexity of the computer model lies in the high-dimensional input space.
$d=53$ random input variables are considered: physical laws essentially, but also initial conditions, material properties and geometrical modeling.
Their probability distributions are either normal or log-normal, and both are truncated.
Such a number of input variables is rather large for the metamodel fitting problem. 
This difficult fit (due to the high dimensionality and small learning sample size) can be made thanks to
the  algorithm of \cite{marioo07}, specifically devoted to this situation.
The obtained Gp metamodel (\ref{Gpmodel}) contains a linear regression part (including $7$ input variables) and a centered Gp model with a generalized exponential 
covariance function (including $6$ input variables).
The reference quality of this Gp model is measured via an additional $1000$-size test sample, which gives $Q_2^{\mbox{\scriptsize ref}} = 0.66$.

\begin{figure}[!ht]
\includegraphics[width=3.5in]{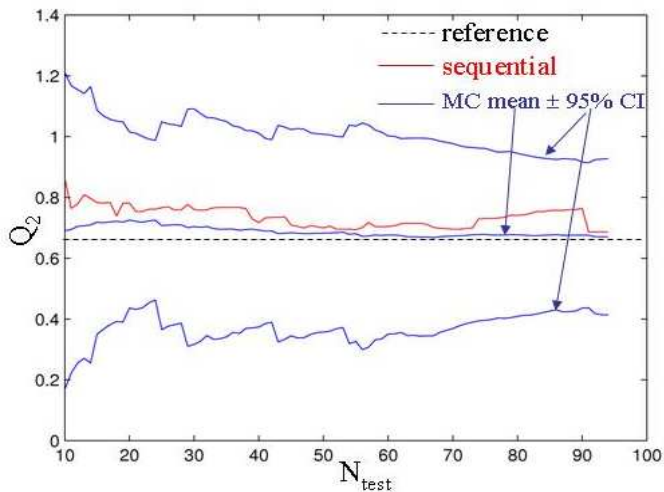}
\caption{For the nuclear safety computer code application, estimation of the metamodel (Gp) predictivity coefficient ($Q_2$) in function of the test sample size $N_{\mbox{\scriptsize test}}$, for two types of validation design: sequential (red) and crude Monte Carlo (blue).}\label{Q2bemuse}
\end{figure}

Figure \ref{Q2bemuse} shows the evolution of the estimated $Q_2$ for test bases with different sizes, ranging from $N_{\mbox{\scriptsize test}}=10$ to $N_{\mbox{\scriptsize test}}=95$.
The sequential validation design gives coarse estimations for all the test design sizes and begins to give precise results for $N_{\mbox{\scriptsize test}} \geq 40$.
Some inadequacies, which remain when $N_{\mbox{\scriptsize test}} \in [75,90]$, have to be finely analyzed in a further work.
In any cases, sequential validation design estimations are clearly less hazardous than using a crude Monte Carlo test sample to validate the metamodel: the $90$\%-confidence intervals obtained using Monte carlo samples show extremely large variation ranges (because of the high dimensionality of the input space: $d=53$).
$Q_2$ estimation using a Monte Carlo test sample can lead to a strongly erroneous result.
Same results have been obtained using optimized LHS for the test design instead of a crude Monte Carlo sample.

\section{Conclusion and future works}

In this paper, we have proposed to look at two practical problems when fitting a metamodel to small-size data samples: the initial design and the validation method choices.
These problems are relevant when a cpu time expensive code has to be replaced by a simplified model with negligible cost (I;e., a metamodel).
Such replacement is useful to resolve optimisation, uncertainty propagation or sensitivity analysis issues.

We have first paid attention to the initial input design.
Our numerical tests concentrate on the popular Gp metamodel and on the LHS.
This type of design, developed thirty years ago, is the most widely used in industrial applications.
We have shown that an excellent way to optimize its properties, in the objective of the best metamodel fit, is to minimize some discrepancy measures, especially the wrap-around $L^2$ discrepancy.
An alternative strategy, if possible, would be to use some adaptive designs \cite{mar08}.
For the Gp metamodel, this kind of design is well-adapted due to the availability of the variance expression (the MSE of the metamodel predictor, see Eq. (\ref{eq_variance})). 

Secondly, we have looked at the metamodel validation process and have shown that the test sample approach can provide erroneous results for small sizes of the test sample.
Moreover, the leave one out approach can strongly underestimate the metamodel predictivity for small sizes of the whole database.
We have proposed to use a recent algorithm, called the sequential validation design, which puts prediction points in the unfilled zones of the learning sample design.
Therefore, a minimal number of points is required to obtain a good estimation of the metamodel predictivity.
Our numerical tests on analytical functions and real application cases have shown that the sequential validation design outperforms the classical metamodel validation methods, especially in high dimensional context.
For our analytical functions, the sequential validation design gives precise estimate of the metamodel predictivity with a test sample size $N_{\mbox{\scriptsize test}} \geq 25$, while for our industrial application, the minimal bound is $N_{\mbox{\scriptsize test}} \geq 40$.

Further works are necessary to more deeply study the validation designs (other test functions with different effective dimensionality and complexity).
Moreover, it would be useful to find a criterion to determine when the sequential validation design can be ended.
Finally, the ultimate goal of such studies will be to define a global strategy of allocating simulation points between the metamodel fitting step and the metamodel validation step.

\section*{Acknowledgments}

This work was supported by the ''Simulation'' program managed by CEA/DEN/DISN.
We thank P. Bazin and A. de Crecy from the CEA Grenoble/DER/SSTH who have provided the CATHARE2 application.


\end{document}